\documentstyle[12pt]{article}
\def\Empty{}
\catcode`\@=11
\def\section{\@startsection {section}{1}{\z@}{-3.5ex plus -1ex minus
-.2ex}{2.3ex plus .2ex}{\large\bf}}

\def\eqalign#1{\,\vcenter{\openup\jot\m@th
  \ialign{\strut\hfil$\displaystyle{##}$&$\displaystyle{{}##}$\hfil
        \crcr#1\crcr}}\,}

\def\mydesc{\list{}{\labelwidth\z@ \itemindent-\leftmargin
\listparindent 1.5em
\let\makelabel\descriptionlabel}}

\def\fnum@figure{{\small Figure \thefigure}}
\def\fakefigure{\def\@captype{figure}}
\long\def\@makecaption#1#2{
    \vskip 10pt
    \def\FCap{#2} \def\NoCap{\ignorespaces}
    \ifx \FCap\NoCap
       \setbox\@tempboxa\hbox{#1}  % This is to avoid the damn colon.
      \else
       \setbox\@tempboxa\hbox{#1: \small \it #2}
    \fi
    \ifdim \wd\@tempboxa >\hsize   % IF longer than one line:
        \unhbox\@tempboxa\par      %   THEN set as ordinary paragraph.
      \else                        %   ELSE  center.
        \hbox to\hsize{\hfil\box\@tempboxa\hfil}
    \fi}

\def\@oddhead{\hbox{}\rightmark \hfil \rm\thepage}% Right heading.
\def\sectionmark#1{\markright {\sc{\ifnum \c@secnumdepth >\z@
      \S\thesection.\hskip 1em\relax \fi #1}}}

\catcode`\@=12

\def\oplabel#1{
  \def\OpArg{#1} \ifx \OpArg\Empty {} \else
        \label{#1}
  \fi}
\def\MakeStEnv#1{
  \newenvironment{#1}[2]{
  \begin{#1St} \oplabel{##1}%
  \global\def\CrntSt{\thetheoremSt}%
  {\rm ##2}%
}{
  \end{#1St} }
}
\MakeStEnv{theorem}
\MakeStEnv{corollary}
\MakeStEnv{proposition}
\MakeStEnv{lemma}
\MakeStEnv{fact}
\MakeStEnv{conjecture}
\MakeStEnv{comment}
\MakeStEnv{remark}
\newenvironment{proof}[1]{
  \def\PfArg{#1}
  \ifx\PfArg\Empty
        \edef\PfArg{\CrntSt}  \fi
 \startproof{\PfArg}%
}{
  \finishproof{\PfArg}
}
\newcommand{\startproof}[1]{
  \medbreak\mbox{}
  {\it Proof of #1:}%
}
\newcommand{\finishproof}[1]{
  \def\FPArg{#1}
  \ifx\FPArg\Empty
        \def\FPArg{\CrntSt}  \fi
  \smallbreak\noindent\makebox[\textwidth]{\hfill\fbox{\FPArg}}
  \medbreak\noindent
}
\newcommand{\hyperbolic}{{\bf H}}

\def\H{\hyperbolic}

\def\til{\widetilde}
\def\hat{\widehat}

\title{The visual core of a hyperbolic $3$-manifold}
\author{James W. Anderson and Richard D. Canary\thanks{Research
supported in part by grants from the National Science Foundation}}
\date{\small{Faculty of Mathematical Studies, University of
Southampton\\ Highfield, Southampton, SO17 1BJ, England}\\
\small{Department of Mathematics, University of Michigan, Ann Arbor,
MI 48109}}

\begin{document}

\maketitle

\section{Introduction}
\label{introduction}

In this note we introduce the notion of the visual core of a hyperbolic
3-manifold $N=\H^3/\Gamma$. One may think of the visual core
as a harmonic analysis analogue of the convex core. Explicitly,
the visual core ${\cal V}(N)$
is the projection to $N$ of all the points in $\H^3$ at which no component
of the domain of discontinuity  of $\Gamma$ 
has visual (equivalently harmonic) measure
greater than half that of the entire sphere at infinity.

We investigate circumstances under which the visual core ${\cal
V}(N')$ of a cover $N'=\H^3/\Gamma'$ of $N$ embeds in $N$, via the
usual covering map $\pi: N'\rightarrow N$. We begin by showing that
the interior of ${\cal V}(N')$ embeds in $N$ when $\Gamma'$ is a
precisely 
QF-embedded subgroup of $\Gamma$, while the visual core itself embeds 
when $\Gamma'$ is a nicely QF-embedded subgroup. We define the notions
of precisely and nicely QF-embedded subgroups of a Kleinian group and
prove these embedding theorems in Section \ref{coverings}.

Applying the results from \cite{ac-cores}, we are able to conclude
that if the algebraic limit of a sequence of
isomorphic Kleinian groups is a generalized web group, then the visual
core of the algebraic limit manifold embeds in the geometric limit
manifold.  This result is part of our ongoing investigation of the
relationship between algebraic and geometric limits of sequences of
isomorphic Kleinian groups. 

\medskip
\noindent
{\bf Theorem \ref{vis-core}: } {\em Let $G$ be a finitely generated,
torsion-free, non-abelian group, let $\{\rho_j\}\subset{\cal D}(G)$ be
a sequence converging algebraically to $\rho\in {\cal D}(G)$, and 
suppose that $\{\rho_j(G)\}$ converges geometrically to
$\widehat\Gamma$. If $\rho(G)$ is a generalized web group, then the
visual core of $N ={\bf H}^3/\rho(G)$ embeds in $\widehat N ={\bf
H}^3/\widehat\Gamma$ under the covering map $\pi:N\rightarrow\widehat
N$.}

\medskip

There are two ways to view Theorem \ref{vis-core}.  On the one hand, one may
think of it as a geometric analogue of the main result from \cite{ac-cores},
which asserts that under the same hypotheses, there is a compact core for the
algebraic limit manifold which embeds in the geometric limit manifold.
On the other hand, Theorem \ref{vis-core} can be thought of as
a  generalization of the result, proven in \cite{accs}, that when the
algebraic limit is a maximal cusp, the convex core of the algebraic
limit manifold embeds in the geometric limit manifold.  In fact, 
when $\Gamma$ is a maximal cusp, the visual and convex cores
of ${\bf H}^3/\Gamma$ coincide.

In Section \ref{combo}, we discuss the relationship between the
visual core and Klein-Maskit combination along component subgroups.
Klein-Maskit combination gives a geometric realization of the
topological operation of gluing hyperbolizable 3-manifolds together
along incompressible surfaces in their boundaries. While the topology
underlying Klein-Maskit combination is well-understood, the geometry
is more mysterious.  For example, the convex core of a summand of a
Klein-Maskit combination need not embed in the resulting
manifold. However, we show in Theorem \ref{KMvismain} that the
(interior of the) visual core of a summand does embed in the resulting
manifold.

There is a relationship between these two investigations, of
limits of sequences of Kleinian groups and of Klein-Maskit
combination, since in the case that the algebraic limit is a
generalized web group, it is shown in \cite{ac-cores} that the
algebraic limit is a summand of a Klein-Maskit decomposition of the
geometric limit.

This paper was completed while the first author was visiting Rice
University, and he would like to thank the department there for their
hospitality.

\section{The visual core}
\label{visual}

Before describing the basic properties of the visual core, we give
some definitions.  A {\em Kleinian group} is a discrete subgroup
of ${\rm PSL}_2({\bf C})$, which we view as acting either on
hyperbolic $3$-space ${\bf H}^3$ as isometries or on the Riemann
sphere $\widehat{\bf C}$ as M\"obius transformations.  The action of
$\Gamma$ partitions $\widehat{\bf C}$ into the {\em domain of
discontinuity} $\Omega(\Gamma)$, which is the largest open subset of
$\widehat{\bf C}$ on which $\Gamma$ acts properly discontinuously,
and its complement the {\em limit set} $\Lambda(\Gamma)$.  The
{\em stabilizer} 
\[ {\rm st}_\Gamma(\Delta)=\{ \gamma \in
\Gamma\: |\: \gamma(\Delta)=\Delta\} \]
of a connected component $\Delta$ of $\Omega(\Gamma)$ is called a {\em
component subgroup} of $\Gamma$.

\medskip

Given a measurable set $X\subset\widehat{{\bf C}}$, consider the
harmonic function $h_X$ on ${\bf H}^3$ defined by setting $h_X(y)$ to
be the proportion of the geodesic rays emanating from $y\in{\bf H}^3$
which end in $X$.  Though we do not explicitly use this formulation,
analytically  we can write $h_X$ in the ball model of hyperbolic
$3$-space as  
\[ h_X(y) =\frac{1}{4\pi}\int_{X}\left(\frac{1-|y|^{2}}
{|y-\zeta|^{2}}\right)^{2} dm(\zeta).\] 
The {\em visual hull} of a Kleinian group $\Gamma$ is then defined to
be 
\[ \widetilde{\cal V}(\Gamma) =\left\{ y\in{\bf H}^3\: |\: h_\Delta
(y)\leq \frac{1}{2}\:\: {\rm for}\:\: {\rm all}\:\: {\rm
components}\:\: \Delta \:\: {\rm of}\:\:\Omega(\Gamma)\right\}.\] 
The {\em visual core} ${\cal V}(N)$ of $N ={\bf H}^3/\Gamma$ is the
quotient $\widetilde{\cal V}(\Gamma)/\Gamma$.  Although ${\cal V}(N)$
is a closed subset of $N$, there is no reason, in general, to suppose
that ${\cal V}(N)$ is a submanifold (or suborbifold) of $N$.

Our first observation is that the visual core of a hyperbolic
3-manifold with finitely generated fundamental group is non-empty
unless its domain of discontinuity is connected and non-empty.

\begin{proposition}{non-empty}{}{}
Let $\Gamma$ be a finitely generated Kleinian group.  Then $\til{\cal
V}(\Gamma)$ is empty if and only if $\Omega(\Gamma)$ is connected and
non-empty.  
\end{proposition}

\begin{proof}{Proposition \ref{non-empty}}  
By definition, $\til{\cal V}(\Gamma)=\H^3$ if and only if the domain
of discontinuity $\Omega(\Gamma)$ of $\Gamma$ is empty. 

If $\Omega(\Gamma)$ is connected and non-empty,
then $\Gamma$ is a function group,
which is a finitely generated Kleinian group whose domain of
discontinuity contains a 
component invariant under the action of the group.  Soma \cite{soma}
shows that $\Gamma$ is then topologically tame, that is 
the orbifold $\H^3/\Gamma$ has a finite manifold cover which is
homeomorphic to the interior of a compact 3-manifold.  Corollary 1 from
\cite{canary-ends} then implies that $\Lambda(\Gamma)$ has measure zero,
so that $h_{\Omega(\Gamma)} (x)=1$ for all $x\in\H^3$.  In particular,
we have that $\til{\cal V}(\Gamma)$ is empty. 

If $\Omega(\Gamma)$ contains at least two components, let $\ell$ be a
geodesic in $\H^3$ whose endpoints at infinity lie in distinct
components $\Delta_1$ and $\Delta_2$ of $\Omega(\Gamma)$. The function
$h_{\Delta_1}$ varies continuously between 0 and 1 on $\ell$, and so
there exists a point $x$ on $\ell$ such that
$h_{\Delta_1}(x)={1\over 2}$.  It follows that $x\in\til{\cal V}(\Gamma)$,
and so $\til{\cal V}(\Gamma)$ is non-empty. 
\end{proof}

It is natural to contrast the definition of the visual core with that
of the convex core.  Recall that the {\em convex hull} $\til{\cal
C}(\Gamma)$ of $\Lambda(\Gamma)$ is obtained from $\H^3$ by removing
each closed hyperbolic half-space which intersects the sphere at
infinity in a closed disk contained in $\Omega(\Gamma)$.  The {\em
convex core} ${\cal C}(N)$ of $N ={\bf H}^3/\Gamma$ is the quotient
$\widetilde{\cal C}(\Gamma)/\Gamma$.  Equivalently, the convex core of
$N$ is the smallest convex submanifold of $N$ whose inclusion is a
homotopy equivalence.  (See Epstein and Marden \cite{epstein-marden}
for further discussion of the convex core.)

The following proposition describes the basic relationship between the
visual and convex cores of a hyperbolic 3-manifold $N$.

\begin{proposition}{two-cores}{}{} Let $N ={\bf H}^3/\Gamma$ be a
hyperbolic $3$-manifold.  Then, its visual core ${\cal V}(N)$ is
contained in its convex core ${\cal C}(N)$. Moreover, the visual core
is equal to the convex core if and only if the boundary
$\partial{\cal C}(N)$ of the convex core is totally geodesic. 
\end{proposition}

\begin{proof}{Proposition \ref{two-cores}} For each point $x$ of ${\bf
H}^3 - \widetilde{\cal C}(\Gamma)$, there exists a hyperplane $H$ in
${\bf H}^3$ containing $x$ so that the circle at infinity $C$ of $H$
bounds a closed disk contained entirely in a component $\Delta$ of
$\Omega(\Gamma)$.  Thus, $h_\Delta (x) > \frac{1}{2}$, which implies
that $x\notin \widetilde{\cal V}(\Gamma)$. Therefore, $\widetilde{\cal
V}(\Gamma)\subset \widetilde{\cal C}(\Gamma)$, which in turn implies
that ${\cal V}(N)\subset{\cal C}(N)$. 

For each point $x$ of $\partial \til{\cal C}(\Gamma)$, there exists a
hyperplane $H$ in ${\bf H}^3$ containing $x$ so that the circle at
infinity $C$ of $H$ bounds an open disk $D$ contained entirely in a
component $\Delta$ of $\Omega(\Gamma)$. If $D$ does not equal
$\Delta$, then $h_\Delta (x) > \frac{1}{2}$, which implies that
$x\notin \widetilde{\cal V}(\Gamma)$.  Therefore, $\til{\cal
C}(N)=\til{\cal V}(N)$ if and only if each component of
$\Omega(\Gamma)$ is a circular disc, which is equivalent to requiring
that $\partial{\cal C}(N)$ be totally geodesic. 
\end{proof}

\section{The visual core and coverings}
\label{coverings}

In this section, we develop a criterion, expressed in terms of limit
sets, which guarantees that the visual core of a cover of a
hyperbolic manifold embeds under the covering map.  This criterion
involves the introduction of two closely related notions of how a
subgroup $\Gamma'$ of a Kleinian group $\Gamma$ sits inside $\Gamma$.

We begin by observing that if $\Gamma'$ is a precisely QF-embedded
subgroup of a Kleinian group $\Gamma$, then the interior of
${\cal V}(N')$ embeds in $N$ under the covering map $\pi: N \rightarrow N'$
(where $N=\H^3/\Gamma$ and $N' ={\bf H}^3/\Gamma'$).  Here, a subgroup
$\Gamma'$ of a Kleinian group $\Gamma$ is {\em precisely QF-embedded}
if, for each $\gamma\in\Gamma -\Gamma'$, there is a component $\Delta$
of $\Omega(\Gamma')$ so that $\gamma(\Lambda(\Gamma'))$ is contained
in $\overline{\Delta}$, ${\rm st}_{\Gamma'}(\Delta)$ is
quasifuchsian, and $\Delta$ is a Jordan domain.

Recall that a {\em quasifuchsian} group is a finitely generated
Kleinian group whose limit set is a Jordan curve and which stabilizes
both components of its domain of discontinuity.
In particular, if $\Delta$ is a component of the domain of
discontinuity $\Omega(\Gamma)$ of a Kleinian group $\Gamma$, 
${\rm st}_{\Gamma'}(\Delta)$ is quasifuchsian, and $\Delta$ is a Jordan domain,
then $\Lambda({\rm st}_{\Gamma'}(\Delta))=\partial \Delta$.
If $\Gamma$ is finitely generated, then a component $\Delta$ of
$\Omega(\Gamma)$ is a Jordan domain if and only if
${\rm st}_{\Gamma}(\Delta)$ is quasifuchsian,  see Lemma 2 of Ahlfors
\cite{ahlfors-structure} and Theorem 2 of Maskit \cite{maskit-boundaries}. 
Hence, a finitely generated subgroup $\Gamma'$ of a Kleinian group
$\Gamma$ is precisely QF-embedded if, for each $\gamma\in \Gamma-\Gamma'$,
there is a component $\Delta$ of $\Omega(\Gamma')$ so that
$\gamma(\Lambda(\Gamma'))$ is contained in $\overline{\Delta}$ and
$\Delta$ is a Jordan domain.  Precisely QF-embedded subgroups arise
naturally in Klein-Maskit combination theory, as we see in Section \ref{combo}.

\begin{proposition}{pecase}{}{} Let $\Gamma$ be a Kleinian group and
let $\Gamma'$ be a precisely QF-embedded subgroup
of $\Gamma$.  Let $N = {\bf H}^3/\Gamma$ and $N' ={\bf H}^3/\Gamma'$,
and let $\pi:N'\rightarrow N$ be the covering map.  Then, $\pi$ is an
embedding restricted to the interior of the visual core ${\cal V}(N')$
of $N'$.
\end{proposition}

\begin{proof}{Proposition \ref{pecase}}
Since the interior of ${\cal V}(N')$ is an open submanifold (or
possibly an open sub-orbifold, in the case that $\Gamma'$ contains
torsion) of $N'$, it suffices to show that $\pi$ is injective on
the interior of ${\cal V}(N')$.  As
${\cal V}(N')$ is covered by $\til{\cal V}(\Gamma')\subset{\bf H}^3$,
it suffices to show that if $\gamma \in \Gamma-\Gamma'$, then
$\gamma({\rm int}(\til{\cal V}(\Gamma'))) \cap {\rm int}(\til{\cal V}(\Gamma'))$
is empty.
 
Let $\gamma$ be any element of $\Gamma-\Gamma'$.  Since $\Gamma'$ is
precisely QF-embedded in $\Gamma$, there exists a component $\Delta_1$
of $\Omega(\Gamma')$ so that
$\gamma(\Lambda(\Gamma'))\subset \overline{\Delta_1}$,
${\rm st}_{\Gamma'}(\Delta_1)$ is quasifuchsian, and $\Delta_1$ is a Jordan
domain.  Since
$U=\widehat{\bf C}-\overline{\Delta_1}$ is a Jordan domain contained
entirely in $\gamma(\Omega(\Gamma'))=\Omega(\gamma\Gamma'\gamma^{-1})$, there
exists a component $\Delta_1'$ of $\gamma(\Omega(\Gamma'))$ such that
$U\subset \Delta_1'$.  In particular,
$\overline{\Delta_1}\cup\Delta'_1 =\widehat {\bf C}$.  Since
$\partial \Delta_1$ is the limit set of the quasifuchsian group
${\rm st}_{\Gamma'}(\Delta_1)$, it has measure
zero, and so $h_{\Delta_1}(x) +h_{\Delta'_1} (x)\ge 1$ for all
$x\in {\bf H}^3$.

Since $h_{\Delta_1}$ is harmonic  and non-constant, it cannot be
locally constant. Thus, if $x$ is in ${\rm int}(\til{\cal V}(\Gamma'))$,
we see that $h_{\Delta_1}(x)< {1\over 2}$.  Since
$h_{\Delta_1'}(x)\ge 1-h_{\Delta_1}(x)>{1\over 2}$, we see that $x$
does not lie in $\gamma({\rm int}(\widetilde{\cal 
V}(\Gamma')))$. Therefore $\gamma({\rm int}(\widetilde{\cal
V}(\Gamma'))) \cap {\rm int}(\widetilde{\cal V}(\Gamma'))$ is empty,
as desired. 
\end{proof}

We next observe that if $\Gamma'$ is nicely QF-embedded, then the
visual core ${\cal V}(N')$ embeds in $N$ under the covering map.
Here, a subgroup $\Gamma'$ of a Kleinian group $\Gamma$ is {\em nicely
QF-embedded} if, for each $\gamma\in\Gamma -\Gamma'$, there is a
component $\Delta$ of $\Omega(\Gamma')$ so that
$\gamma(\Lambda(\Gamma'))$ is contained in $\overline{\Delta}$, 
${\rm st}_{\Gamma'}(\Delta)$ is quasifuchsian, $\Delta$ is a Jordan
domain, and $\gamma(\Lambda(\Gamma'))\cap\partial\Delta\ne
\partial\Delta$. 
More simply, if $\Gamma'$ is a finitely generated subgroup of
a Kleinian group $\Gamma$, then $\Gamma'$ is nicely QF-embedded
if, for each $\gamma\in\Gamma -\Gamma'$, there is a
component $\Delta$ of $\Omega(\Gamma')$ so that
$\gamma(\Lambda(\Gamma'))$ is contained in $\overline{\Delta}$,
$\gamma(\Lambda(\Gamma'))\cap\partial\Delta\ne \partial\Delta$,
and $\Delta$ is a Jordan domain.  A nicely QF-embedded subgroup is always
precisely QF-embedded, though the converse need not hold. Nicely embedded
QF-subgroups occur naturally in the study of algebraic and geometric
limits, as we see in Section \ref{limits}.

\begin{proposition}{necase}{}{} Let $\Gamma$ be a Kleinian group and
let $\Gamma'$ be a finitely generated, nicely
QF-embedded subgroup of $\Gamma$.  Let $N = {\bf H}^3/\Gamma$ and $N'
={\bf H}^3/\Gamma'$, and let $\pi:N'\rightarrow N$ be the covering
map. Then, $\pi$ is an embedding restricted to the visual
core ${\cal V}(N')$ of $N'$.
\end{proposition}

\begin{proof}{Proposition \ref{necase}}
We argue much as  in the proof of Proposition \ref{pecase} to show
that $\pi$ is injective on ${\cal V}(N)$.  Let $\gamma$ be any element
of $\Gamma-\Gamma'$.  Since $\Gamma'$ is nicely QF-embedded in
$\Gamma$, there exists a component $\Delta_1$ of $\Omega(\Gamma')$
which is a Jordan domain, so that $\gamma(\Lambda(\Gamma'))\subset
\overline{\Delta_1}$ and $\partial \Delta_1-\gamma(\Lambda(\Gamma'))$
is non-empty.  
Thus $\widehat{\bf C}-\overline{\Delta_1}$ is a proper subset of some
component $\Delta_1'$ of $\gamma(\Omega(\Gamma'))$.  In particular,
$\overline{\Delta_1}\cup\Delta'_1 =\widehat{\bf C}$ and
$\Delta_1\cap\Delta'_1\ne\emptyset$.  Since $\partial \Delta_1$ is the
limit set of a quasifuchsian group, it has measure zero. Thus,
$h_{\Delta_1}(x)+h_{\Delta'_1}(x) >1$ for all $x\in {\bf H}^3$.

If $x\in\widetilde{\cal V}(\Gamma')$, then $h_{\Delta_1}(x)\le {1\over 2}$.
So, we see that $h_{\Delta_1'}(x) > 1-h_{\Delta_1}(x)\ge {1\over 2}$,
which implies that $x$ does not lie in $\gamma(\til{\cal V}(\Gamma'))$.
Thus, $\gamma(\widetilde{\cal V}(\Gamma')) \cap \til{\cal V}(\Gamma')$
is empty, which proves that $\pi$ is
injective on the visual core ${\cal V}(N')$. 

To verify that $\pi$ is an embedding restricted to ${\cal V}(N')$, it
only remains to check that $\pi$ is proper.  If not, then there must 
exist a sequence $\{ x_j\}$ of points in $\partial {\cal V}(N')$ which
exits every compact subset of $N'$, but such that $\{\pi(x_j)\}$
converges to a point $x$ in $N$. By passing to a subsequence, we may
assume that $d(x_j,x_{j+1})\ge 1$ and $d(\pi(x_j),x)\le {1\over 3j}$
for all $j$.  

Let $\{ \tilde{x}_j\}$ be a sequence
of lifts of $\{ x_j\}$ to ${\bf H}^3$. Since
$d(\pi(x_j),\pi(x_{j+1}))<{1\over j}$ 
and $d(x_j,x_{j+1})\ge 1$, for each $j$ there exists an element
$\gamma_j\in \Gamma-\Gamma'$ such that
$d(\til x_j,\gamma_j(\til x_{j+1}))< {1\over j}$.
Since $\Gamma'$ is nicely QF-embedded, there exists,
for each $j$, a component $\Delta_j$ of $\Omega(\Gamma')$ which is
a Jordan domain whose closure contains $\gamma_j(\Lambda(\Gamma'))$.
Since $\til x_{j+1}\in \til{\cal V}(\Gamma)$,
$\gamma_j(\til{\cal V}(\Gamma)) \cap \til{\cal V}(\Gamma) =\emptyset$, and
$\gamma_j(\Lambda(\Gamma))\subset \Delta_j$, we see that
$h_{\Delta_j}(\gamma_j(\til x_{j+1}))>{1\over 2}$.
Since $\til x_j\in \til{\cal V}(\Gamma')$,
$h_{\Delta_j}(\til x_j)\le {1\over 2}$. So, by continuity, there exists
a point $\til q_j$ between $\til x_j$ and $\gamma_j(\til x_{j+1})$
such that $d(\til q_j,\til x_j)< {1\over j}$ and 
$h_{\Delta_j}(\til q_j)={1 \over 2}$. Thus, 
$\{ q_j\}$ is a sequence of points in $\partial {\cal V}(N')$ which
exits every compact subset of $N'$, but such that $\{\pi(q_j)\}$
converges to a point $x$ in $N$.

Since  $\Gamma'$ is finitely generated, there exist only
finitely many inequivalent components of $\Omega(\Gamma')$, so we may assume
(by choosing different lifts, passing to a subsequence, and
relabelling) that there exists a fixed component $\Delta_0$ of
$\Omega(\Gamma')$ which is a Jordan domain, so that
$h_{\Delta_0}(\tilde{q}_j) =\frac{1}{2}$ for all $j$. 
Let $\Gamma_0={\rm st}_{\Gamma'}(\Delta_0)$, $N_0 ={\bf H}^3/\Gamma_0$
and $p:{\bf H}^3\rightarrow N_0$ be the covering map.
Since $h_{\Delta_0}(\tilde{q}_j)={1\over 2}$, we conclude that
$\tilde{q}_j\in \widetilde{{\cal V}}(\Gamma_0)$.
Let $y_j =p(\tilde{q}_j)$.  Since $\{ q_j\}$ exits
every compact subset of $N'$, $\{ y_j\}$ must exit every compact subset 
of $N_0$.  Proposition \ref{two-cores} guarantees that the sequence
$\{ y_j\}$ lies entirely in the convex core $C(N_0)$ of $N_0$.
Since $\Gamma'$ is finitely generated and $\Delta_0$ is a Jordan domain,
$\Gamma_0$ is quasifuchsian.  Therefore, the $\epsilon$-thick part
of the convex core,
$$C(N_0)_\epsilon =  \{ y\in C(N_0)|\ {\rm inj}_{N_0}(y)\ge \epsilon\},$$
is compact for all $\epsilon>0$, see Bowditch \cite{bowditch}.
(Here, ${\rm inj}_{N_0}(y)$ denotes the injectivity radius of
the point $y$ in $N_0$.) Thus, 
${\rm inj}_{N_0}(y_j) \rightarrow 0$, which implies that
${\rm inj}_{N}(\pi(x_j)) \rightarrow 0$. However, this contradicts the
fact that $\{ \pi(x_j)\}$ converges in $N$.
\end{proof}

\noindent{\bf Remarks:}
(1) One may think of Proposition \ref{pecase} as an analogue of
Proposition 6.1 of \cite{ac-cores}, which asserts that if $\Gamma'$ is
a finitely generated, torsion-free, precisely embedded generalized web
subgroup of $\Gamma$, then there is a compact core for $N'$ which
embeds (via the covering map $\pi:N\rightarrow N'$) in $N$. That
result may be generalized, using the same techniques as in
\cite{ac-cores}, to show that if $\Gamma'$ is a finitely generated,
torsion-free precisely  QF-embedded subgroup of $\Gamma$, then there
is a compact core for $N'$ which embeds (via the covering map $\pi$)
in $N$. This generalization is the more direct topological analogue of
Proposition \ref{pecase}.

(2) The arguments in the proofs of Propositions
\ref{pecase} and \ref{necase} may be used to show that larger subsets
embed. Let $\til{\cal W}(\Gamma')$ be the set of points $x\in\H^3$ such that
$h_\Delta(x)\le {1\over 2}$ for every component $\Delta$ of
$\Omega(\Gamma')$ such that ${\rm st}_{\Gamma'}(\Delta)$ is quasifuchsian,
and set ${\cal W}(N')=\til{\cal W}(\Gamma')/\Gamma'$. Then one can
adapt the proof of Proposition \ref{pecase} to show
that if $\Gamma'$ is a precisely QF-embedded subgroup of $\Gamma$,
then the interior of ${\cal W}(N')$ embeds in $N$. Similarly, one can
adapt the proof of Proposition \ref{necase} to show that if $\Gamma'$
is a finitely generated, nicely QF-embedded subgroup of $\Gamma$,
then ${\cal W}(N')$ embeds in $N$.

(3) Note that the definitions for a precisely QF-embedded and of a
nicely QF-embedded subgroup $\Gamma'$ of a Kleinian group $\Gamma$
both make sense for an 
infinitely generated subgroup $\Gamma'$.  In fact, Proposition
\ref{pecase} as stated holds for infinitely generated, precisely
QF-embedded subgroups.  The reason that in the definitions we require
both that the component $\Delta$ of $\Omega(\Gamma')$ be a Jordan
domain and that ${\rm st}_{\Gamma'}(\Delta)$ be quasifuchsian is that
it is possible, by taking the Klein combination of a quasifuchsian
group with an infinitely generated Kleinian group with trivial
component subgroups, to construct an infinitely generated Kleinian
group whose component subgroups are all quasifuchsian but the
components of its domain of discontinuity are not all simply
connected, and so in particular cannot all be Jordan domains. 

\section{Algebraic and geometric limits}
\label{limits}

In an earlier paper \cite{ac-cores}, we proved that if the algebraic limit
of a sequence of isomorphic Kleinian groups is a generalized web group,
then it is a nicely QF-embedded subgroup of the geometric limit. In that
paper, we used this result to establish that there is a compact core for the
algebraic limit manifold which embeds in the geometric limit manifold,
thus obtaining ``topological'' information about how the algebraic
limit sits inside the geometric limit.  In this section, we use the
results of the previous section to obtain ``geometric'' information
about how the algebraic limit sits inside the geometric limit.

We briefly recall the basic definitions from the theory of algebraic
and geometric limits.  We refer the interested reader to J\o
rgensen and Marden \cite{jorgensen-marden} for more details. Given a
finitely generated group $G$, let ${\cal D}(G)$ denote the space of
discrete, faithful representations of $G$ into ${\rm PSL}_2({\bf
C})$.  A sequence $\{\rho_i\}$ in ${\cal D}(G)$ converges {\em
algebraically} to $\rho$ if $\{\rho_i(g)\}$ converges to $\rho(g)$ for
each $g\in G$.  

A sequence $\{\Gamma_j\}$ of Kleinian groups converges {\em
geometrically} to a Kleinian group $\widehat\Gamma$ if every element
of $\widehat\Gamma$ is the limit of a sequence
$\{\gamma_j\in\Gamma_j\}$ and if every accumulation point of every
sequence $\{\gamma_j\in\Gamma_j\}$ lies in $\widehat\Gamma$.  If $G$ is
not virtually abelian and if $\{\rho_i\}$ converges to $\rho$ in
${\cal D}(G)$, then there is a subsequence $\{\rho_j(G)\}$ of
$\{\rho_i(G)\}$ which converges geometrically to a Kleinian group
$\hat\Gamma$ which contains $\rho(G)$.

In this note, we restrict ourselves to sequences $\{\rho_n\}$ in
${\cal D}(G)$ so that $\{\rho_n\}$ converges algebraically to some
$\rho\in {\cal D}(G)$ and so that $\{\rho_n(G)\}$ converges
geometrically to $\hat\Gamma$.  The Kleinian group $\rho(G)$ is the
{\em algebraic limit} of $\{\rho_n\}$, and $\hat\Gamma$ is the {\em
geometric limit} of $\{\rho_n(G)\}$.  If $G$ is torsion-free, we refer
to ${\bf H}^3/\rho(G)$ as the {\em algebraic limit manifold}
and to ${\bf H}^3/\hat\Gamma$ as the {\em geometrical limit manifold}.  Since
$\rho(G)\subset\hat\Gamma$, there is a natural covering map $\pi: {\bf
H}^3/\rho(\Gamma)\rightarrow {\bf H}^3/\hat\Gamma$.  In order to
understand the relationship between the algebraic and geometric limit,
it is important to understand how $\rho(G)$ ``sits inside'' $\hat\Gamma$,
which is closely related to understanding the covering map $\pi$.

A finitely generated Kleinian group $\Gamma$ is called a {\em generalized web
group} if $\Omega(\Gamma)$ is non-empty and if every component
subgroup of $\Gamma$ is quasifuchsian (or equivalently, if every
component of  $\Omega(\Gamma)$ is a Jordan domain). Theorem A from
\cite{ac-cores} asserts that if the algebraic limit is a generalized
web group, then it is a nicely QF-embedded subgroup of the geometric limit. 

\begin{theorem}{nicely-embedded}{(Theorem A of \cite{ac-cores})}{} Let
$G$ be a finitely generated, torsion-free, non-abelian group, let
$\{\rho_j\}$ be a sequence in ${\cal D}(G)$ converging algebraically
to $\rho\in {\cal D}(G)$, and suppose that $\{\rho_j(G)\}$ converges
geometrically to $\widehat\Gamma$.  If $\rho(G)$ is a generalized web
group, then $\rho(G)$ is a nicely QF-embedded subgroup of
$\widehat\Gamma$.
\end{theorem}

One may combine Theorem \ref{nicely-embedded} with Proposition \ref{necase}
to  obtain ``geometric'' information about how the algebraic limit sits
within the geometric limit in this case.

\begin{theorem}{vis-core}{}{} Let $G$ be a finitely generated,
torsion-free, non-abelian group, let $\{\rho_j\}\subset{\cal D}(G)$ be
a sequence converging algebraically to $\rho\in {\cal D}(G)$, and
suppose that $\{\rho_j(G)\}$ converges geometrically to
$\widehat\Gamma$. If $\rho(G)$ is a generalized web group, then the
visual core of $N ={\bf H}^3/\rho(G)$ embeds in $\widehat N ={\bf
H}^3/\widehat\Gamma$ under the covering map $\pi:N\rightarrow\widehat
N$.
\end{theorem}

\noindent
{\bf Remarks:} (1) One may think of Theorem \ref{vis-core} as one way
to generalize Proposition 3.2 from \cite{accs}, which shows that if the
algebraic limit is a maximal cusp, then the convex core of the
algebraic limit manifold embeds in the geometric limit manifold under
the covering map. In fact, one may view our Theorem \ref{vis-core}
and the result from \cite{ac-cores} that asserts that, under the
same assumptions, a compact core for the algebraic limit manifold
embeds in the geometric limit manifold, as two different generalizations
of Proposition 3.2 from \cite{accs}.  

(2) In general, even if the algebraic limit is a
generalized web group, the convex core of the algebraic limit manifold
need not embed in the geometric limit manifold.

(3) The examples given in \cite{ac-books} illustrate the point that
the visual core of the algebraic limit manifold need not embed in the
geometric limit manifold in the case that the algebraic limit is not a
generalized web group.

\section{Klein-Maskit Combination}
\label{combo}

In this section, we discuss the relationship between the visual core
and the operation of Klein-Maskit combination.  We restrict our entire
discussion to Klein-Maskit combination along component subgroups.
For a more complete discussion of Klein-Maskit combination
see Maskit \cite{maskit-book}.  In this setting,
we see that the interior of the visual core of a summand of a
Klein-Maskit decomposition of a hyperbolic 3-manifold embeds in the
manifold.

There are two types of Klein-Maskit combination.  The first, type I,
corresponds topologically to gluing 2 hyperbolic 3-manifolds together
along incompressible components of their conformal boundary. The
second, type II, corresponds topologically to gluing together two
incompressible components of the conformal boundary of a single
hyperbolic 3-manifold.

The following theorem summarizes the relevant properties of Klein-Maskit
combination of type I along a component subgroup (see Theorem VII.C.2
in \cite{maskit-book}).

\begin{theorem}{combination_I}{}{(Klein-Maskit combination I)}
Let $\Gamma_1$ and $\Gamma_2$ be Kleinian groups, and let
\newline
$\Phi =\Gamma_1\cap\Gamma_2$.  Suppose that $\Phi$ is a quasifuchsian
group which is a component subgroup of both $\Gamma_1$ and $\Gamma_2$,
and that $\Lambda(\Gamma_1)$ and $\Lambda(\Gamma_2)$ lie in  the
closures of different components of $\Omega(\Phi)$.
Then,
\begin{enumerate}
\item
$\Gamma =\langle\Gamma_1,\Gamma_2\rangle$ is a Kleinian group isomorphic
to the amalgamated free product of $\Gamma_1$ and $\Gamma_2$ along
$\Phi$; 
\item $\Gamma_1$ and $\Gamma_2$ are nicely QF-embedded subgroups of
$\Gamma$; 
\item If $\gamma \in \Gamma-\Gamma_i$, then
$\gamma(\Lambda(\Gamma_i))$ is contained in a component $\Delta$ of
$\Omega(\Gamma_i)$ which is $\Gamma_i$-equivalent to the component
$\Delta_i$ of $\Omega(\Gamma_i)$ bounded by $\Lambda(\Phi)$. Moreover,
$\partial \Delta-\gamma(\Lambda(\Gamma_i))$ is non-empty; and
\item
$\H^3/\Gamma$ is homeomorphic to the manifold (or orbifold) obtained  from
$(\H^3\cup \Delta_1)/\Gamma_1$ and $(\H^3\cup \Delta_2)/\Gamma_2$ by 
identifying $\Delta_1/\Phi$ with $\Delta_2/\Phi$.
\end{enumerate}
\end{theorem}

In this case, we say that $\Gamma_1$ is a {\em summand of a simple type
I Klein-Maskit decomposition} of $\Gamma$.  Combining property (2)
of Theorem \ref{combination_I} with Proposition \ref{necase} yields
the following result: 

\begin{proposition}{KMvis1}{}{}
Let $\Gamma_1$ be a finitely generated Kleinian group which is a
summand of a simple type I Klein-Maskit decomposition of $\Gamma$, and
set $N_1=\H^3/\Gamma_1$ and $N=\H^3/\Gamma$.  Then, the visual core
${\cal V}(N_1)$ of $N_1$ embeds in $N$ (via the covering map
$\pi:N_1\to N$). 
\end{proposition}

If $\Gamma_1$ is not finitely generated, the techniques in the proof
of Proposition \ref{necase} may be adapted to show that
${\cal V}(N_1)$ still embeds in $N$. More simply, one may combine
property (2) of Theorem \ref{combination_I} with Proposition
\ref{pecase} to obtain the following weaker result:
 
\begin{proposition}{KMvis1gen}{}{}
Let $\Gamma_1$ be a summand of a simple type I Klein-Maskit decomposition
of  a Kleinian group $\Gamma$, and
set $N_1=\H^3/\Gamma_1$ and $N=\H^3/\Gamma$.  Then, the interior of
the visual core ${\cal V}(N_1)$ of $N_1$ embeds in $N$
(via the covering map $\pi:N_1\to N$).
\end{proposition}

Moreover, if $\Gamma$ is a simple type I
Klein-Maskit combination of $\Gamma_1$ and $\Gamma_2$,
we may find a larger subset of $N_1$ which embeds in $N$.  Consider the sets
\[ \til{\cal X}_i= \left\{ x\in {\bf H}^3\: | \: h_\Delta (x) \le
{1\over 2} \: {\rm for\ all\ components}\ \Delta \ {\rm of} \
\Omega(\Gamma') \ {\rm equivalent \ in}\ \Gamma_i \  {\rm to} \
\Delta_i \right\}. \]
Clearly, the visual hull $\til{\cal V}(\Gamma_i)$
of $\Gamma_i$ is contained  in $\til{\cal X}_i.$
Let ${\cal X}_i=\til{\cal X}_i/\Gamma_i\subset
N_i={\H^3/\Gamma_i}$.  Condition (3) above and the arguments in the
proof of Proposition \ref{necase} can then be used to show
that ${\cal X}_i$ embeds in $N=\H^3/\Gamma$.

\bigskip

The following theorem summarizes the relevant properties of Klein-Maskit
combination of type II along a component subgroup
(see Theorem VII.E.5 in \cite{maskit-book}).

\begin{theorem}{combination_II}{}{(Klein-Maskit combination II)} Let
$\Gamma_1$ be a Kleinian group, and let $\Delta$ and $\Delta'$ be
components of $\Omega(\Gamma_1)$ which are Jordan domains. Suppose
that $\Phi={\rm st}_{\Gamma_1}(\Delta)$ and $\Phi'={\rm st}_{\Gamma_1}(\Delta')$
are quasifuchsian and are not conjugate by an element of $\Gamma_1$.
Let $\gamma$ be a M\"{o}bius transformation which conjugates
$\Phi'$ to $\Phi$, and assume that $\gamma(\Lambda(\Gamma_1))$
and $\Lambda(\Gamma_1)$ lie in the closures of different components
of $\Omega(\Phi)$. Then,
\begin{enumerate}
\item
$\Gamma =\langle\Gamma_1,\gamma\rangle$ is a Kleinian group isomorphic
to the HNN-extension of $\Gamma_1$ with stable letter $\gamma$
and associated subgroups $\Phi$ and $\Phi'$;
\item $\Gamma_1$ is a precisely QF-embedded subgroup of $\Gamma$;
\item If $\gamma \in \Gamma-\Gamma_1$, then $\gamma(\Lambda(\Gamma_1))$
is contained in a component of $\Omega(\Gamma_1)$ which is
$\Gamma_1$-equivalent to either $\Delta$ or $\Delta'$;
\item
$\H^3/\Gamma$ is homeomorphic to the manifold  (or orbifold) obtained  from
$(\H^3\cup \Delta \cup \Delta')/\Gamma_1$ by identifying
$\Delta/\Phi$ with $\Delta'/\Phi'$ by the homeomorphism determined by
$\gamma$. 
\end{enumerate}
\end{theorem}

In this case, we say that $\Gamma_1$ is a {\em summand of simple type
II Klein-Maskit decomposition} of $\Gamma$. Combining property (2)
of Theorem \ref{combination_II} with Proposition \ref{pecase} yields
the following result: 
 
\begin{proposition}{KMvis2}{}{}
Let $\Gamma_1$ be a summand of simple type II KIein-Maskit
decomposition of $\Gamma$, and set $N_1=\H^3/\Gamma_1$ and
$N=\H^3/\Gamma$.  Then, the interior of the visual core ${\cal
V}(N_1)$ of $N_1$ embeds in $N$ (via the covering map
$\pi:N_1\to N$).
\end{proposition}
 
As in the type I situation, if $\Gamma_1$ is a summand of
simple type II Klein-Maskit decomposition of $\Gamma$, we may find
a larger subset of $N_1$ which embeds in $N$. Consider the set 
\[ \til{\cal X}= \left\{ x\in {\bf H}^3\: | \: h_\Delta (x) \le
{1\over 2} \: {\rm for\ all\ components}\ \Delta \ {\rm of} \
\Omega(\Gamma') \ {\rm equivalent \ in}\ \Gamma_1\  {\rm to} \ \Delta\
{\rm or}\ \Delta'\right\}, \]
which clearly contains the visual hull $\widetilde{\cal V}(\Gamma_1)$
of $\Gamma_1$.  Let ${\cal X}=\til{\cal X}/\Gamma_1\subset N_1$.
Condition (3) above and the arguments in the proof of Proposition \ref{pecase}
can be used to show that the interior of ${\cal X}$ embeds in $N=\H^3/\Gamma$. 

\bigskip

In general, if a Kleinian group $\Gamma$ can be built from $\Gamma_1$
and a collection of other Kleinian groups by repeatedly performing
Klein-Maskit combinations of types I and/or II along component
subgroups, we say that $\Gamma_1$ is a {\em summand of a Klein-Maskit
decomposition} of $\Gamma$. By applying Propositions \ref{KMvis1gen} and
\ref{KMvis2}, we obtain the following summation of the results 
of this section.

\begin{theorem}{KMvismain}{}{}
Let $\Gamma_1$ be a summand of a Klein-Maskit decomposition of
$\Gamma$. If $N_1=\H^3/\Gamma_1$ and $N=\H^3/\Gamma$, then the
interior of the visual core ${\cal V}(N_1)$ of $N_1$ embeds in $N$
(via the covering map $\pi:N_1\to N$).
\end{theorem}

\noindent
{\bf Remarks:} (1) The definition of the visual core was suggested
by Thurston's reproof of the Klein-Maskit combination theorems. In
our notation, Thurston shows that in the type I decomposition,
$N=\H^3/\Gamma$ is obtained from ${\cal X}_1$ and ${\cal X}_2$ by
identifying points in their boundaries. (In general, one must be a
little careful since ${\cal X}_i$ need not be a submanifold.)
Similarly, in the type II situation he shows that $N$ is obtained from
${\cal X}$ by identifying points in the boundary. (This proof is
discussed in outline in Section 8 of Morgan \cite{Morgan}.) 

(2) Notice that if $\Gamma_1$ is a summand of a simple type II
Klein-Maskit decomposition of a Kleinian group $\Gamma$, then
$\Gamma_1$ is a precisely QF-embedded subgroup of $\Gamma$, but is not
a nicely QF-embedded subgroup.  In this same case, the interior
of the visual cover of $N_1$ embeds in $N$, but the visual core itself
does not.

(3) Corollary D of \cite{ac-cores} asserts that if the algebraic limit
of a sequence of isomorphic Kleinian groups is a generalized web
group, then it is a summand of a Klein-Maskit decomposition of the
geometric limit. Hence, there is a close relationship between 
Theorems \ref{vis-core} and \ref{KMvismain}.

\end{document}